\documentclass[11pt]{amsart}
\usepackage{setspace}
\usepackage{graphicx}
\usepackage{amsmath}
\usepackage{amsfonts}
\usepackage{amssymb}
\begin{document}

\centerline{\LARGE \bf Blaschke, Osgood, Wiener, Hadamard and the Early}

\smallskip
\centerline{\LARGE \bf Development of Modern Mathematics in China}

\bigskip
\centerline{\Large Chuanming Zong}

\vspace{0.6cm}
\centerline{\begin{minipage}{12.8cm}
{\bf Abstract.} In ancient times, China made great contributions to world civilization and in particular to mathematics. However, modern sciences including mathematics came to China rather too late. The first Chinese university was founded in 1895. The first mathematics department in China was formally opened at the university only in 1913. At the beginning of the twentieth century, some Chinese went to Europe, the United States of America and Japan for higher education in modern mathematics and returned to China as the pioneer generation. They created mathematics departments at the Chinese universities and sowed the seeds of modern mathematics in China. In 1930s, when a dozen of Chinese universities already had mathematics departments, several leading mathematicians from Europe and USA visited China, including Wilhelm Blaschke, George D. Birkhoff, William F. Osgood, Norbert Wiener and Jacques Hadamard. Their visits not only had profound impact on the mathematical development in China, but also became social events sometimes. This paper tells the history of their visits.
\end{minipage}}

\vspace{0.4cm}

\vspace{1cm}
\noindent
{\Large\bf 1. Introduction and Background}

\vspace{0.3cm}
\noindent
In 1895, Peiyang University was founded in Tianjin as a result of the westernization movement. In 1898, Peking University was founded in Beijing. They are the earliest universities in China.

On May 28, 1900, Britain, USA, France, Germany, Russia, Japan, Italy and Austria invaded China to suppress the Boxer Rebellion\footnote{At the end of the nineteenth century, after the second Anglo-Chinese war, foreign missionaries were under particular protection in China. Gradually, their followers also benefited from the protection. In 1897, there was a conflict between the followers of a church and the local villagers in Shandong Province. When the villagers thought that they were treated unfairly, they asked a famous boxer Sanduo Zhao to administer justice. They got the expected justice! Afterwards, organized boxer gangs for the purpose of fighting with the missionaries and their followers were quickly formed and spread from Shandong to the neighbor provinces, finally to Tianjin and Beijing. Many churches were destroyed and more than 100,000 people from both sides were killed. At that time, Empress Dowager Ci Xi (1835-1908) was in power. Since the rebels took a slogan {\it supporting the Qing rulers and opposing foreign powers}, they obtained tacit permission from Ci Xi. In June 1900, more than 100,000 boxer rebels gathered in Beijing. When thousands of missionaries and their followers fled to the foreign embassies, the boxer rebels sieged the embassy area in Beijing.} which was an organized movement against foreign missionaries and their followers. On August 14, 1900, more than 20,000 soldiers of the allied force occupied the Chinese capital Beijing and the emperor's family with some high ranking officials and servants escaped to Xi'an. On August 8, 1901, when the Qing government agreed all the requirements of the eight nations including paying them huge indemnity (the Boxer Indemnity), the allied force withdrew from Beijing.

On December 28, 1908, after four years of hard diplomatic negotiation, American government finally admitted that the indemnity was overpaid and decided to return the overpayment of the Boxer Indemnity to China. That was 11.96 million US dollars. Both sides agreed to set up an education foundation with this money to support Chinese students to study in America. According to the agreement, in the first four years from 1909, China could send one hundred students each year and thereafter fifty students a year until the fund was exhausted. Meanwhile, Tsinghua School was founded in 1909 to take charge of the selection and preparation of the students.

The education foundation set up by the returned Boxer Indemnity was crucial for the Chinese modern mathematics. It led to the first Chinese Ph.D Mingfu Hu (1891-1927) in Harvard in 1917, educated the first generation of the great Chinese mathematical educators such as Li-Fu Chiang (Chan-Chan Tsoo, 1890-1978) and King-Lai Hiong (1893-1969), and educated the first generation of globally known Chinese mathematicians such as Loo-Keng Hua (1910-1985), Shiing-Shen Chern (1911-2004) and Pao Lu Hs$\ddot{\rm u}$ (1910-1970).

In 1912, the first mathematical department in China was set up at Peking University with two professors: Zuxun Feng (1880-1940) and Junji Hu (1886-??). The next year it enrolled two students as a start. Professor Feng taught analytical courses such as calculus, function theory and differential equations, and professor Hu taught the others such as Euclidean geometry and integer theory.

Zuxun Feng enrolled in Peking University as a undergraduate student in 1902. Two years later, he was chosen to study mathematics in Japan supported by Peking University. He graduated from Kyoto Imperial University in 1909, specialized in differential equations. This made him the first Chinese who ever had a complete university education in mathematics. Junji Hu went to Japan in 1903 with his brother. He studied Mining Metallurgy and graduated from Tokyo Imperial University. At that time, he learnt some modern mathematics. In particular, when he took the teaching position at Peking University, he visited Tokyo Imperial University once again to improve his mathematics.

In 1917, two more professors joined the mathematical faculty: Fen Qin (1887-1973) and Renfu Wang (1886-1959). Qin obtained a Master degree in astronomy in 1909 and Wang obtained a Bachelor degree in mathematics in 1913, both from Harvard. At that time, the math department had seven faculty members and more than twenty students. It offered a dozen of mathematical courses such as Euclidean geometry, calculus, function theory, differential equations, harmonic functions, theoretic physics, abstract algebra, modern geometry, group theory and number theory.

In 1919, Nankai University was founded in Tianjin, China. The next year it set up a mathematical department and appointed Li-Fu Chiang, who obtained his Ph.D from Harvard under the supervision of J. L. Goolidge, as the chair professor. It was the second mathematical department in China. In fact he was the only faculty member of the math department for four years. He alone taught various courses in analysis, geometry and algebra. During that hard time, he educated several promising students such as Chin-Nien Liu (1904-1968), Tsai-Han Kiang (1902-1994) and Yu-Cheng Shen (1901-1978). Afterwards, all Liu, Kiang and Shen obtained their Ph.D from Harvard and made distinguished contribution to modern mathematical development in China.

In 1920, National Southeast University was founded in Nanjing, China. In the next year it set up a mathematical department and appointed King-Lai Hiong (1893-1969), who obtained a Master degree from Universit$\acute{\rm e}$ de Montpellier in 1919, as its founding chair professor. He designed the courses, prepared lecture notes and taught analytic geometry, spherical geometry, calculus, analytic functions, differential geometry and differential equations. In 1922, Zi-Xie Duan (1890-1969) joined the mathematical faculty as a professor. Then the department could offer more mathematical courses. Duan went to France in 1913. He obtained a Master degree in mathematics from Lyon in 1920.

In 1927, Tsinghua University set up a mathematical department. King-Lai Hiong and Zhi-Fan Zheng (1887-1963) were appointed as the first two professors, with professor Zheng as chairman. Zhi-Fan Zheng studied mathematics at Cornell and graduated in 1910. Afterwards, he visited Harvard to improve his mathematics for one year.

In three years, two more professors joined the faculty: Dan Sun (Guang-Yuan Sun, 1900-1979) and Ko-Chuen Yang (1896-1973). Both Dan Sun and Ko-Chuen Yang obtained their Ph.D from Chicago in 1928, the first wrote a thesis in differential geometry under the supervision of E. P. Lane and the second wrote his thesis in number theory under the supervision of L. E. Dickson. Then, the mathematical faculty at Tsinghua became one of the strongest math faculties in China.

In 1930s, China had already more than forty universities and several competitive mathematical departments such as the departments at Peking University, Tsinghua University, Zhejiang University, Central University, Nankai University, Chiao Tung University, Kwang Hua University and Wuhan University (see Zhang \cite{Zhang-99}). Up to 1930, eighteen Chinese obtained their Ph.D in mathematics, as listed in the following table. They played crucial roles for introducing modern mathematics to China (see Zhang \cite{Zhang-16}).

\bigskip
\centerline{\begin{tabular}{|c|c|c|c|}
\hline
Name & University & Time & Field \\
\hline
Mingfu Tah Hu & Harvard & 1917 & Analysis \\
\hline
Chan-Chan Tsoo & Harvard & 1919 & Geometry \\
\hline
Jung Sun & Syracuse & 1921 & Algebra \\
\hline
David Yule & Harvard & 1922 & Mathematical Logic \\
\hline
Bing-Chin Wong & UC Berkeley & 1922 & Geometry \\
\hline
Shih-Luan Wei & G\"ottingen & 1925 & Analysis \\
\hline
Zhao-An Zeng & Columbia & 1925 & Geometry \\
\hline
Kun-Ching Chu & G\"ottingen & 1927 & Analysis \\
\hline
Ko-Chuen Yang & Chicago & 1928 & Number Theory \\
\hline
Dan Sun & Chicago & 1928 & Geometry \\
\hline
Chin-Yi Chao & Lyon & 1928 & Analysis \\
\hline
Wei-Kwok Fan & Lyon & 1929 & Analysis \\
\hline
Kien-Kwong Chen & Tohoku & 1929 & Analysis \\
\hline
Tsun-Shien Lian & Lyon & 1930 & Analysis \\
\hline
Tsai-Han Kiang & Harvard & 1930 & Topology \\
\hline
Chin-Nien Liu & Harvard & 1930 & Analysis \\
\hline
Hung-Chi Chang & Michigan & 1930 & PDE \\
\hline
Shu-Ting Liu & Michigan & 1930 & Analysis \\
\hline
\end{tabular}}

\bigskip
In 1918, Dr. Mingfu Tah Hu published his thesis at {\it Transactions of the American Mathematical Society}. It was the first mathematical research paper published by Chinese. Up to 1930, about one hundred mathematical papers were published by Chinese authors.

\vspace{0.8cm}
\noindent
{\Large \bf 2. Painlev$\acute{\rm\bf e}$ and Russell's Visits to China}

\vspace{0.3cm}\noindent
The first western mathematicians to visit China were Paul Painlev$\acute{\rm e}$ (1863-1933) and Emile Borel (1871-1956) from France, though their visit was not aimed at mathematics. Painlev$\acute{\rm e}$ had twice been the Prime Minister of France, in 1917 and in 1925. In 1920, he made an official visit to China, leading a French delegation including Emile Borel. At that time, Borel was president of Ecole Normale Superieure.

On July 1, 1920, Painlev$\acute{\rm e}$ was awarded an honorary doctorate by Peking University in Beijing. On that occasion, he gave a talk entitled {\it Mathematical Progress}, which was the first mathematical talk given by a foreign speaker in China. Unfortunately, no mathematical talk was given by E. Borel. In Shanghai, Painlev$\acute{\rm e}$ presented a talk {\it Science and Education in China} at the Chinese Society of Sciences, in which he appealed Chinese scholars to organize societies in their own specialized fields. His appeal stimulated Chinese mathematicians to establish the Chinese Mathematical Society in 1935 (see Li \cite{li-14}).

On 12 October 1920, Bertrand Russell (1872-1970) arrived in Shanghai with his girlfriend. During his first three weeks in China, Russell delivered several public lectures in Shanghai, Hangzhou, Wuhan and Changsha on topics ranging from Einstein's relativity theory to education and social problems. His lectures had influence on many Chinese social activists at that time and stimulated debates at newspapers. Zedong Mao (1893-1976) attended his talk in Changsha. Early in November, Russell arrived at his ultimate destination, Peking University in Beijing. According to the schedule, he should lecture on problems of philosophy for six months.

At the beginning of March 1921, invited jointly by the Society for Mathematics and Physics at Peking University and by the Society for Mathematics, Physics and Chemistry at Beijing Normal University, Bertrand Russell agreed to give four lectures on mathematical logic. The first lecture was given on 8 March to an audience of about 150 professors and college students. The second lecture was announced for 15 March. Unfortunately, the day before, Russell came down with pneumonia and almost died. For nearly six weeks he was confined to his bed. He eventually recovered, but was still very weak when he left Beijing on 10 July 1921. Because of his sudden illness, Russell's subsequent lectures on mathematical logic were all cancelled. Nevertheless, his visit and talks stimulated several Chinese to pursue mathematical logic (see Xu \cite{Xu-03}).

\vspace{0.8cm}
\noindent
{\Large \bf 3. Blaschke and Sperner at Peking University}

\vspace{0.3cm}\noindent
In April 1932, as a part of his round-the-world mathematical tour to India, China, Japan and the United States of America, the German mathematician Wilhelm Blaschke (1885-1962) visited Peking University for two weeks. His book {\it Reden und Reisen eines Geometers} contains an account of the whole tour. During the visit, he gave a series of talks on differential geometry and integral geometry.

At that time, there were very few people in Peking (Beijing) area who could understand differential geometry. One of them was professor Dan Sun (Guang-Yuan Sun, 1900-1979) from Tsinghua University, who obtained his Ph.D from Chicago in 1928 with a thesis in differential geometry under the supervision of E. P. Lane. Another two were Shiing-Shen Chern and Da-Ren Wu (1908-1997) who were graduate students of professor Sun working on their Master degrees at Tsinghua.

Blaschke's visit was not as sensational as Russell's. Nevertheless, it was a big event for the Chinese mathematical community. According to \cite{Blaschke-61}, within the two weeks he was invited to more than a dozen of banquets, from which he gained much weight. A lot of people came from other cities to attend his talks, even though many of them could not understand a word. On this occasion Chern first met Blaschke. He attended all the lectures, took detailed notes, and made his decision to study in Hamburg.

After his brief visit, Blaschke recommended E. Sperner as a visiting faculty member at Peking University for two years. He taught geometric foundation and topology, ran advanced seminars, and supervised graduate students. At that time, the mathematical department had about seventy undergraduates and four graduate students. During his visit, Dr. Sperner also gave several talks at Tsinghua University and Nankai University.

When Shiing-Shen Chern finished his Master degree, he succeeded in the selection examination to study in abroad financed by the Boxer Indemnity. Usually, the students supported by the Boxer Indemnity should pursue their study and research in the United States. However, attracted by Blaschke, Chern decided to make his Ph.D in Humburg. Surprisingly, his application was approved by the Boxer Indemnity authority. In 1934, Shiing-Shen Chern arrived in Hamburg and enrolled as a Ph.D student of Blaschke. From there, he gradually stepped up the world mathematical stage.

In fact, Da-Ren Wu also went to Hamburg a couple of years later. It was him who first introduced integral geometry into China.

\vspace{0.8cm}
\noindent
{\Large \bf 4. Birkhoff and Osgood at Peking University}

\vspace{0.3cm}\noindent
Modern mathematics was introduced into the United States of America in 1870s. The American Mathematical Society was founded in 1888. Both George D. Birkhoff (1884-1944) and William F. Osgood (1864-1943) belonged to the first generation of the globally known American mathematicians.

In 1919, Li-Fu Chiang obtained his Ph.D from Harvard under the supervision of J. L. Goolidge. He was the second Chinese who ever made a Ph.D in mathematics. When Chiang studied at Harvard, he was once an assistant of professor William F. Osgood who was chairman of the mathematics department at that time. In 1919, Dr. Chiang was appointed professor and founding chairman of the mathematical department of Nankai University.

At Nankai University, Tsai-Han Kiang (1902-1994) was one of the earliest students majoring in mathematics. Supported by a Boxer Indemnity fellowship, he went to USA and obtained his Ph.D from Harvard in 1930 under the supervision of M. Morse. In 1931, Dr. Kiang was appointed a professor at Peking University. Three years later, he was appointed head of the mathematics department.

In April 1934, arranged by Tsai-Han Kiang, professor George D. Birkhoff from Harvard visited Peking University. During his visit, he gave a series of talks on several solutions in quantum mechanics, differential equations of dynamics, four-colour problem, and aesthetic measurement. His talks were so fashionable that they attracted great interest among both professors and students alike. However, since they were too advanced, almost no audience could understand them.

In 1933, William F. Osgood retired from Harvard University, after being a professor there for thirty years. Recommended by Li-Fu Chiang and invited by Tsai-Han Kiang, he came to China and taught as a visiting professor at Peking University from 1934 to 1936. During his visit, contrary to Birkhoff, professor Osgood taught various basic courses on mechanics, real functions and complex functions. As assistants of the professor, Pao Lu Hs$\ddot{\rm u}$ and Shu-Ben Sun made careful notes of the lectures. Afterwards, his lecture notes {\it Functions of Real Variables} and {\it Functions of a Complex Variable} were published by Peking University (see Ding, Yuan and Zhang \cite{Ding-93}).

The long term visits of Sperner, Birkhoff and Osgood, in particular their basic mathematics courses, made great inspiration and encouragement to the Chinese colleagues and students. Their lecture notes were rare and valuable text books for years.

\vspace{0.8cm}
\noindent
{\Large \bf 5. Wiener at Tsinghua University}

\vspace{0.3cm}\noindent
From 1924 to 1930, Yuk-Wing Lee studied electrical engineering at Massachusetts Institute of Technology. In 1930, he obtained his Ph.D there under the supervision of professor V. Bush. During his MIT time, with the recommendation of his supervisor, Lee became an assistant and close friend of Norbert Wiener (1894-1964). In 1932, Dr. Lee returned to China and was appointed professor at the electrical engineering department of Tsinghua University. At that time, professor King-Lai Hiong was the head of the mathematics department. In 1934, based on the recommendation of professor Lee and professor Hiong, the president of Tsinghua University invited Norbert Wiener for a long term visit. At the age of only forty, professor Wiener already made a name as a mathematician. He was elected to the National Academy of Sciences, USA, in 1933.

The Wiener family arrived at Tsinghua on 15 August, 1935. With the help of Yuk-Wing Lee, his family set up in the campus of Tsinghua University. In the next two semesters, professor Wiener taught several courses on Fourier series, Fourier integrals, and Lebesgue integrals. To prepare the courses, based on his suggestions, the university bought the related reference books in advance. His classes were very successful and attracted many faculty members and students alike, not only from Tsinghua, but also from other universities in Beijing area. In particular, his classes paid more attention to raising problems and inspiring ideas rather than simply deducing theorems.

Loo-Keng Hua never had a university education. However, through self-study he became one of the most important modern mathematicians in China. When Wiener visited Tisinghua University, Hua was an assistant at the mathematical faculty. Of course, he attended all the courses of Wiener, kept discussing mathematics with him, and became friends with him. Gradually, by reading his papers, Wiener was impressed by his mathematical talent and offered to recommend him to G. H. Hardy in Cambridge. In 1936, supported by a Boxer Indemnity fellowship, Loo-Keng Hua went to Cambridge University for two years, which became a turning point for his mathematical career.

Besides his lectures in mathematics department, professor Wiener continued to work with Dr. Lee on problems of electric-circuit design. They tried follow in the footsteps of Bush in making an analogy-computing machine, but to gear it to the high speed of electrical circuits instead of to the much lower one of mechanical shafts and integrators. In fact, when Wiener accepted the invitation, the chairman of the electrical engineering department asked him for help to buy an analogy-computing machine in USA for Tsinghua University. Unfortunately, it was too expensive and Tsinghua could not afford it. Anyway, Wiener spoke highly of his visit to China. He wrote in his book \cite{Wiener-56} \lq\lq If I were to take any specific boundary point in my career as a journeyman in science and as in some degree an independent master of the craft, I should pick out 1935, the year of my China trip, as that point." In the book, he also wrote in detail about his impression of the Chinese people, their culture, beliefs and daily life.

The Wiener family left China on 19 May, 1936, for the International Congress of Mathematicians in Oslo. He maintained lasting friendships with his Tsinghua colleagues, in particular with Loo-Keng Hua and Yuk-Wing Lee, and provided his support whenever it was needed. Based on his recommendation, John von Neumann was much interested to visit China. Unfortunately, the Japanese invasion in 1937 made it no longer possible.

\vspace{0.8cm}
\noindent
{\Large \bf 6. Hadamard at Tsinghua University}

\vspace{0.3cm}\noindent
On 22 March, 1936, Jacques Hadamard (1865-1963) and his wife arrived in Shanghai by passenger ship the {\it Queen of Asia}. When the ship stopped at the wharf, representatives of the Chinese Mathematical Society, the Chinese Physical Society and the Sino-French Friendship Association got on board
to greet them. In the evening, president Yuan-Pei Cai of Academia Sinica gave a banquet to welcome them. Many famous scientists and scholars in Shanghai attended the welcome banquet. In the following days, professor Hadamard gave a couple of public lectures at Chiao Tung University and Sino-French Friendship Association, respectively.

At that time, both Hadamard and his wife were about seventy. Therefore, before their departure, they were a little worried about their health
and the security in Beijing. To make sure the security situation, he wrote to Norbert Wiener who had been there already for one semester. Of course, he obtained a warm confirmation.

On 7 April, 1936, they arrived in Beijing by train. The Hadamards were greeted by the president of Tsinghua University, the dean of science and
the chairman of mathematics department at the railway station. At that time, professor King-Lai Hiong was the chairman who motivated the invitation and made the arrangement.

According to the plane, professor Hadamard would give twenty lectures on definite problems of partial differential equations. It was the first time that partial differential equations was systematically taught in China. Of course, he also gave public mathematical talks on different occasions. For example, on the 25th anniversary of Tsinghua University, he lectured on {\it some reflections on the role of mathematics.} His talks were great inspiration to the audience, in particular to the young faculty members and students.

Everybody knows that the prime number theorem was first proved by Jacques Hadamard and Charles-Jean de la Vall$\acute{\rm e}$e Poussin respectively in 1896. Although his mathematical interest already shifted to partial differential equations and other subjects, Hadamard himself was a great master in number theory. When he was at Tsinghua, he quickly got to know Loo-Keng Hua and impressed by his mathematical talent and his persistency in the subject. While discussing the Waring problem with Hua, Hadamard introduced Vinogradov's work to him. Afterward, Hua wrote to Vinogradov and made a life long friendship with him. Along with Norbert Wiener, professor Hadamard also persuaded chairman Hiong to support Loo-Keng Hua go abroad to improve his mathematics.

Besides Loo-Keng Hua, professor Hadamard also helped several other young Chinese. With his recommendations, both Xin-Mou Wu and Chi-Tai Chuang went to Paris supported by the Boxer Indemnity. Wu studied partial differential equations with professors H. Villat and Hadamard. He returned to China in 1951 and became one of the PDE pioneers in China. Chuang obtained his Ph.D in 1938 under the supervision of professor G. Valiron from University of Paris. He returned to China in 1939 and became a leading expert in complex functions in China.

When Hadamards arrived at Tsinghua, the Wiener family was still there. The two great men had a lot of mathematical discussions and the two families had many happy times together. Both Hadamard's biography \cite{mazya-98} and Wiener's autobiography \cite{Wiener-56} gave accounts on their experiences in Beijing. Wiener wrote \lq\lq We used to go to town to visit the Hadamards and sometimes Margaret and I, or Lee and the two of us, used to go down into the tangled, squalid streets of the socalled Chinese city (as opposed to the rectangular Tatar city) to rummage in the antique shops. There we would often come across ancestor portraits which show dignified Chinese gentlemen or ladies, in stiff poses, with hands on the knees, dressed in marvelous silken gowns, which for the men were robes of office, civil or military. For all their pomp and stiffness, it was common for the faces in these pictures to be of a remarkable fineness, humor, and sensitivity. we found one such ancestor portrait which was so like Professor Hadamard himself, with his somewhat sparse, stringy beard, his hooked nose, and his fine, sensitive features, that it would have been completely adequate to identify him and to pick him out of a large assembly of people. There was, it is true, a very slight slant to the eyes, and a very slight sallowness of complexion, but not enough to confuse the identification. We bought this picture and gave it to its likeness. He appreciated it very much, but I don't think that Mme. Hadamard cared for it."

When the Chinese Mathematical Society was founded in 1935, it decided to launch two journals {\it Acta Mathematica Sinica} and {\it Journal of Mathematics}, the first for original works and the second for introductional reports. To support the Chinese mathematical community, both professor Wiener and professor Hadamard published papers at the early issues of the Acta.

On 25 June, 1936, professor Hadamard finished his Chinese visit and left Beijing for Paris with his wife by the Trans-Siberian Railway. His visit had lasting impact on the mathematical development in China.

\vspace{0.8cm}\noindent
{\bf Acknowledgements.} This work is supported by the National Natural Science Foundation of China (NSFC11921001) and the National Key Research and Development Program of China (2018YFA0704701).

\vspace{0.5cm}
\bibliographystyle{amsplain}

\bigskip\medskip\noindent
Chuanming Zong, Center for Applied Mathematics, Tianjin University, Tianjin, China.

\noindent
Email: cmzong@math.pku.edu.cn

\end{document}